\documentstyle[11pt,twoside]{article}

\catcode`\@=11

\font\twtymy = msym10 \@magscale4
\font\svtnmy = msym10 \@magscale3
\font\frtnmy = msym10 \@magscale2
\font\twlmy  = msym10 \@magscale1
\font\elvmy  = msym10 \@halfmag
\font\tenmy  = msym10
\font\ninmy  = msym9
\font\egtmy  = msym8
\font\sevmy  = msym7
\font\sixmy  = msym6
\font\fivmy  = msym5

\newfam\my

\@addfontinfo\@xxvpt{\textfont\my\twtymy
		\scriptfont\my\twtymy\scriptscriptfont\my\svtnmy}
\@addfontinfo\@xxpt{\textfont\my\twtymy
		\scriptfont\my\frtnmy\scriptscriptfont\my\twlmy}
\@addfontinfo\@xviipt{\textfont\my\svtnmy
		\scriptfont\my\twlmy\scriptscriptfont\my\tenmy}
\@addfontinfo\@xivpt{\textfont\my\frtnmy
		\scriptfont\my\tenmy\scriptscriptfont\my\sevmy}
\@addfontinfo\@xiipt{\textfont\my\twlmy
		\scriptfont\my\egtmy\scriptscriptfont\my\sixmy}
\@addfontinfo\@xipt{\textfont\my\elvmy
		\scriptfont\my\egtmy\scriptscriptfont\my\sixmy}
\@addfontinfo\@xpt{\textfont\my\tenmy
		\scriptfont\my\sevmy\scriptscriptfont\my\fivmy}
\@addfontinfo\@ixpt{\textfont\my\ninmy
		\scriptfont\my\sixmy\scriptscriptfont\my\fivmy}
\@addfontinfo\@viiipt{\textfont\my\egtmy
		\scriptfont\my\sixmy\scriptscriptfont\my\fivmy}
\@addfontinfo\@viipt{\textfont\my\sevmy
		\scriptfont\my\sixmy\scriptscriptfont\my\fivmy}
\@addfontinfo\@vipt{\textfont\my\sixmy
		\scriptfont\my\sixmy\scriptscriptfont\my\fivmy}
\@addfontinfo\@vpt{\textfont\my\fivmy
		\scriptfont\my\fivmy\scriptscriptfont\my\fivmy}

\newcommand{\Bbb}[1]{\fam\my\relax#1}
\def\OO{{\cal O}}
\def\RR{{\Bbb R}}
\def\ZZ{{\Bbb Z}}

\def\CC{{\Bbb C}}

\def\PP{{\Bbb P}}
\def\FF{{\Bbb F}}
\def\CF{{\cal F}}

\def\CL{{\cal L}}

\def\@item[#1]{\if@noparitem \@donoparitem
  \else \if@inlabel \indent \par \fi
         \ifhmode \unskip\unskip \par \fi
         \if@newlist \if@nobreak \@nbitem \else
                        \addpenalty\@beginparpenalty
                        \addvspace{0.5\@topsep} \addvspace{-\parskip}\fi
           \else \addpenalty\@itempenalty \addvspace{-2pt}          \fi
    \global\@inlabeltrue
\fi
\everypar{\global\@minipagefalse\global\@newlistfalse
          \if@inlabel\global\@inlabelfalse \hskip -\parindent \box\@labels
             \penalty\z@ \fi
          \everypar{}}\global\@nobreakfalse
\if@noitemarg \@noitemargfalse \if@nmbrlist \refstepcounter{\@listctr}\fi \fi
\setbox\@tempboxa\hbox{\makelabel{#1}}%
\global\setbox\@labels
 \hbox{\unhbox\@labels \hskip \itemindent
       \hskip -\labelwidth \hskip -\labelsep
       \ifdim \wd\@tempboxa >\labelwidth
                \box\@tempboxa
          \else \hbox to\labelwidth {\makelabel{#1}}\fi
       \hskip \labelsep}\ignorespaces}

\def\ENDBOX{\hfill\mbox{$\Box$}}
\def\qed{\ENDBOX}

\newtheorem{theorem}{Theorem}[section]%
\newtheorem{theorem*}{Theorem}%
\newtheorem{proposition}[theorem]{Proposition}%
\newtheorem{corollary}[theorem]{Corollary}%
\newtheorem{Lemma}[theorem]{Lemma}%
\def\thethoremcounter{\thetheorem}%
\newcommand\Thmcounter{\refstepcounter{theorem}{\bf \thethoremcounter}}%

\newenvironment{Def}%
	{\begin{trivlist}%
	\item[]{\bf Definition\ \Thmcounter\ } \begingroup \rm}%
	{ \endgroup \end{trivlist}\par}
\newenvironment{remark}%
	{\begin{trivlist}%
	\item[]{\it Remark\ \Thmcounter\ } \begingroup \rm}%
	{ \endgroup \end{trivlist}\par}
\newenvironment{example}%
	{\begin{trivlist}%
	\item[]{\it Example\ \Thmcounter\ } \begingroup \rm}%
	{ \endgroup \end{trivlist}\par}
	{\begin{trivlist}%
	\item[]{\bf \Thmcounter\ } \begingroup \rm}%
	{ \endgroup \end{trivlist}\par}

\newenvironment{pf}%
	{\begin{trivlist}%
	\item[]{\it Proof. } \begingroup \rm}%
	{ \ENDBOX \endgroup \end{trivlist}\par}

\catcode`\@=12

\author{A. Steffens}
\date{\today}

\title{On the stability of the tangent bundle of  Fano manifolds
\footnote{
       This article contains results from the author's dissertation
       which was prepared at the graduate program "Complex Manifolds" at
        the university of Bayreuth.
       The author wants to express his gratitude to his advisors Prof. T.
Peternell and
      Prof. M. Schneider.}}

\begin{document}
\maketitle

\section*{Introduction}
A smooth variety $X$ over the field of complex numbers $\CC$ is called
{\em Fano} if its anticanonical divisor $-\!K_X$ is ample.
Stability (in the sense of Mumford and Takemoto)
with respect to $-K_X$ of the tangent bundle $T_X$ can be considered as an
algebraic analogue to the existence of a K\"ahler-Einstein metric on $X$,
since the result of Kobayashi \cite{K} and L\"ubke \cite{Lubke} shows
that the existence of a K\"ahler-Einstein metric implies the stability
of  the tangent bundle.
But the converse is not true, e.g.~ $\PP^2$ blown up in two
points has stable tangent bundle, which do not
admit a K\"ahler-Einstein metric, cf.~ \cite{Ma}.

By Tian's solution of Calabi's conjecture for Del-Pezzo surfaces \cite{Ti}
and by \cite{Fahlaoui} we have a complete picture in dimension 2:
If $X$ is a Del-Pezzo surface, then $X$ has stable tangent bundle $T_X$,
unless $X$ is isomorphic to $\PP^1\times\PP^1$, or $\PP^2$ blown-up in
a point. In both cases the relative tangent bundle $T_{X/\PP^1}$ of a
canonical projection to $\PP^1$ is a destabilising subsheaf of $T_X$.
If the dimension of $X$ is $\geq 3$, then the existence of a K\"ahler-Einstein
 metric remains an open question.

In this article, our main results are as follows:

\begin{theorem*}\label{dim 3 main 1}
Let $X$ be a Fano 3-fold with $b_2\geq 2$. Assume that the tangent bundle $T_X$
 of $X$
is not stable.

Then the relative tangent sheaf $T_{X/Y}$ of
 a contraction  $f:X\longrightarrow Y$ of an extremal face on $X$
 is a destabilising
subsheaf of $T_X$.
\end{theorem*}

\begin{theorem*}\label{dim 3 main 2}
{}From the 87 deformation classes of Fano 3-folds with $b_2\geq 2$,
cf.~ \cite{Mori-Mukai.1,Mori-Mukai.2} the
members of
\begin{center}
\begin{tabular}{cl}
68 &  deformation classes have stable tangent bundle,\\
12 & deformation classes have semistable (but not stable) tangent bundle and\\
7 & deformation classes have unstable tangent bundle.
\end{tabular}
\end{center}
\end{theorem*}
For a detailed description of the deformation classes whose members
have semistable or unstable tangent bundle see
theorem \ref{dim 3 b2 groesser 1} below.

Our main tool in proving theorem \ref{dim 3 main 2} is Mori theory.
By Mori theory we understand the results and techniques
concerning the cone of curves on a manifold $X$ whose canonical divisor $K_X$
is not
numerically effective.

The proofs of theorem 1 and 2 use the classification of Fano-3-folds.
If $\dim X\geq 4$, then the problem of the stability of the tangent bundle
seems hopeless, if one wants to use classification.
But one may expect that theorem 1 holds in any dimension.

\section{Preliminaries}
A smooth connected variety $X$ over the field of complex numbers $\CC$ is
called simply a {\em manifold}. All manifolds are assumed to be projective,
unless otherwise stated. $K_X$ denotes the canonical divisor of a
normal variety $X$.

Assume $X$ smooth and set $n=\dim X$.
Let $H$ be an ample line bundle on $X$. If $\CF$ is a torsion free coherent
sheaf on $X$ we define $\mu(\CF)$ to be $c_1(\CF).H^{n-1}/rk(\CF)$.
We call $\CF$ {\em semistable} (resp. {\em stable}) if for all proper
subsheaves
$\CF^\prime$ of $\CF$ with $0\le rk(\CF^\prime) \le rk(\CF)$
we have $\mu(\CF^\prime)\leq\mu(\CF)$
(resp. $\mu(\CF^\prime)<\mu(\CF)$).

Let $X$ be a normal variety of dimension $n$. We use the following
notation:
\begin{enumerate}
\item[] $N^1(X):=(\{ \mbox{ Cartier divisors on }X\}/\equiv)\otimes\RR$
\item[] $N_1(X):=(\{ \mbox{ 1-cycles on }X\}/\equiv)\otimes\RR$
\item[] $\overline{N\!E}(X):=$ the closure of the convex cone generated
by effective 1-cycles in $N_1(X)$.
\end{enumerate}
Here the symbol $\equiv$ means  numerical equivalence and the
symbol $\sim$ will denoted  linear equivalence.

\begin{Def}
(1) A curve $C$ on $X$ is called {\em extremal} if
\begin{enumerate}
\item[(a)] $(K_X.C)<0$,
\item[(b)] given $u,v\in\overline{N\!E}(X)$ then $u,v\in\RR_+[C]$ if
 $u+v\in\RR_+[C]$.
\end{enumerate}
If $C$ is an extremal curve on $X$, then the set $R=\RR_+[C]$ is called an
{\em extremal ray} on $X$.\\
2) Let $H$ be a nef Cartier divisor on $X$. The set
$F:= H^\perp\cap\overline{N\!E}\setminus\{0\}$ is called an {\em extremal face}
if $F$ is entirely contained in the set
$\{z\in N_1(X)\ |\ \ (K_X.z)<0\ \}$.

\end{Def}

\begin{theorem}[Cone theorem (Mori, Kawamata, Koll\'ar
\protect{\cite{Mori.1,KMM}})] \hfill\\
Assume that $X$ has only canonical singularities.
Fix an ample divisor $\CL$. Then for any $\varepsilon >0$, there exist
extremal curves $\ell_1,\ldots,\ell_r$ such that
$$\overline{N\!E}(X)=\sum_{i=1}^r \RR_+[\ell_i]
+\overline{N\!E_\varepsilon}(X).$$
Here $\overline{N\!E_\varepsilon}(X):=\{\ z\in\overline{N\!E}(X)\ | \ (K_X.z)>
-\varepsilon
(\CL.z)\ \}$.
\end{theorem}

\begin{theorem}[Contraction theorem (Shokurov, Kawamata \cite{KMM})]\hfill\\
Let $F$ be an extremal face of $\overline{N\!E}(X)$. Assume that $X$ has only
 canonical singularities. Then there exists a morphism
$\varphi=\mbox{cont}_F:X\longrightarrow Y$
onto a normal projective variety Y, such that:
For any irreducible curve $C$ on $X$ the image $\varphi(C)$ is a point
if and only if $[C]\in F$.
\end{theorem}

\section{Fano varieties with $b_2=1$}
Fano 3-folds with $b_2=1$ are classified by Iskovskih
\cite{Iskovskih.1,Iskovskih.2}. There are 18 classes of Fano 3-folds with
$b_2=1$
up to deformation.

\begin{remark}
Let $X$ be a Fano manifold of dimension $n$. By a criterion for
stability \cite{Hoppe}, the tangent bundle $T_X$ of $X$ is stable with respect
to $(-K_X)$,
if one of the following equivalent conditions is fulfilled:
\begin{enumerate}
\item[$(A_i)$]
$H^0(X,{\Omega^i \otimes L^{-1}})  =  0$ for
all $L\in Pic(X)$ with  $L.(-K_X)^{n-1} \geq -\frac{i}{n}(-K_X)^n$.
\item[$(B_i)$] $H^0(X,{\bigwedge^i T_X \otimes L^{-1}})  =  0$ for
all $L\in Pic(X)$ with
 $L.(-K_X)^{n-1} \geq \frac{i}{n}(-K_X)^n$.
\end{enumerate}
Stability is granted when all conditions $(A_i)$ or $(B_{n-i})$,
$1\leq i \leq n-1$, hold.
\end{remark}

{}From now on assume that $b_2(X)=1$. Let $L$ be the ample generator of
$Pic(X) \simeq \ZZ$. Then we have that $-K_X=rL$ with $1\leq r\leq n+1$,
where the integer $r$ is called the {\em index} of $X$.
By the Kobayashi-Ochiai characterisation of projective space and
hyperquadrics \cite{KoOc}, we have:
$$
 r =n+1 \ \Leftrightarrow \ X \simeq \PP^n \mbox{ and} \ \
 r = n   \ \Leftrightarrow \ X \simeq Q_n \subset \PP^{n+1}
$$

\begin{remark}
If $X$ is $\PP^n$ or $Q_n$, then on may verify the conditions $(A_i)$
directly.
\end{remark}

\begin{remark}
Let $X$ be a Fano manifold with $b_2=1$ and $L$ the ample generator
of $Pic(X)$.
Then we have:
\begin{enumerate}
\item
If the index $r$ of $X$ is $1$, then the conditions $(A_i)$ are fulfilled
(cf. \cite[Theorem 3]{Reid.2}).
\item
 $H^0(X,{\Omega^1_X \otimes L^{m}})=0$ for $m \leq 0$.
In particular the condition $(A_1)$ is fulfilled in any case.
\item
 If the index $r$ of $X$ is $\leq n$,
then the condition $(A_{n-1})$ is fulfilled.
\end{enumerate}
\end{remark}

\begin{pf}
1),2)
Since $1\!\leq\! i \!\leq\! n-1$, we have  $i\frac{1}{n} < 1$.
If $m < 0$, then we have $H^0(X,{\Omega^i\! \otimes\! L^m})\!=\!0$
by Kodaira-Nakano vanishing theorem.
Since $-K_X$ is ample, it follows by the Kodaira vanishing theorem that
$h^0\!(\Omega^i_X)\!=\!h^i\!(\OO_X)\!=\!h^{n-i}\!(\omega_X)=0$.

\noindent
3)
Since the condition $(A_{n-1})$ is equivalent to $(B_1)$, it suffices to
show that $H^0(X,{T_X \otimes L^{-m}})  =  0$, for $m \geq 1$.
But this is a consequence of \cite[Theorem 1]{Wahl}, because the index
of $X$ is different from $n+1$.
\end{pf}

\begin{corollary}
Let $X$ be a Fano 3-fold with $b_2=1$.
Then the tangent bundle of $X$ is stable.
\end{corollary}

\section{Fano 3-folds with $b_2\geq 2$}
The proof of theorem \ref{dim 3 main 1} and theorem \ref{dim 3 main 2}
is a by-product of the proof of the following:
\begin{theorem}\label{dim 3 b2 groesser 1}
Let $X$ be a Fano 3-fold.\\
i) The members of the deformation classes in the following list have semistable
tangent  bundle.
\begin{enumerate}
\item[(1)]
the blow-up of $\PP^3$ with center a line.
\item[(2)]
the  blow-up of $Y$  with center
two exceptional fibers $\ell$ and $\ell^\prime$ of the blow-up
$\Phi : Y\to \PP^3$ such that $\ell$ and $\ell^\prime$ lie on the same
irreducible component of the exceptional set of $\Phi$. Here
$\Phi : Y\to \PP^3$ is the blow-up of $\PP^3$ with center two disjoint lines
in $\PP^3$.
\item[(3)]
the product of a Del-Pezzo surface (i.e Fano 2-fold) with $\PP^1$.
\end{enumerate}
ii) The members of the deformation classes in the following list have unstable
tangent  bundle.
\begin{enumerate}
\item[(1)]
$V_7$, that is, the $\PP^1$-bundle $\PP(\OO\oplus\OO(1))$ over $\PP^2$.
\item[(2)]
the blow-up of the Veronese cone $W_4\subset\PP^6$ with center the vertex,
that is $\PP(\OO\oplus\OO(2))$ over $\PP^2$.
\item[(3)]
the blow-up of $V_7$ with center a line on the exceptional divisor
$D\simeq\PP^2$ of the blow-up $V_7\to\PP^3$.
\item[(4)]
the blow-up of $V_7$ with center the strict transform of a line passing
through the center  of the blow-up $V_7\to\PP^3$.
\item[(5)]
the blow-up of the cone over a smooth quadric surface in $\PP^3$ with center
the vertex, that is, the $\PP^1$-bundle
$\PP(\OO\oplus\OO(1,1))$ over $\PP^1\times\PP^1$.
\item[(6)]
the blow-up of $\PP^1\times\FF_1$ with center $\{t\}\times e$, where
$t\in\PP^1$
and $e$ is an exceptional curve of the first kind on $\FF_1$.
\item[(7)]
the  blow-up of $\widetilde{\PP}(L)$  with center
two exceptional lines  of the blow-up $\widetilde{\PP}(L)\to\PP^3$.
Here
$\widetilde{\PP}(L)\to \PP^3$ is the blow-up of $\PP^3$ with center a line
in $\PP^3$.
\end{enumerate}
iii)  If $X$ is not contained in a deformation class
listed as above, then $T_X$ is stable with respect to the anticanonical
divisor $-K_X$.
\end{theorem}
\begin{pf}
Instead of presenting here the long proof of theorem \ref{dim 3 b2 groesser 1},
we will treat some special cases and examples.
For the proof of theorem \ref{dim 3 b2 groesser 1}, we will refer the
reader to \cite{St}.
The plan of the proof is as follows:

\noindent
Step 1. Vanishing results for
$H^0(X,T_X\otimes\CL^{-1})$ and $H^0(X,\Omega^1_X\otimes\CL^{-1})$.

\noindent
Step 2. Direct check of stability for the list of
families with $b_2=2$.

\noindent
Step 3.
Reduction of the cases with $b_2\geq 3$ to those studied
at Step 2, or to lower dimensional vanishing statements.
\end{pf}

\noindent
Products of Fano manifolds\\
Let $Y_1, Y_2$ be two Fano manifolds of dimension $n_1$ and $n_2$ respectively.
Then $X\!=\!Y_1\!\times\!Y_2$ is a Fano manifold of dimension $n=n_1 + n_2$.
By an easy computation, one gets
$\mu(T_X)= \mu({\pi_1^*}T_{Y_1}) =  \mu({\pi_2^*}T_{Y_2})$.
It is a well known fact that a vector bundle $E_1\oplus E_2$
is semistable if and only if $E_1$ and  $E_2$ are semistable vector bundles
with $\mu(E_1)=\mu(E_2)$. Thus, we have proved:

\begin{proposition}
Let $Y_1$, $Y_2$ be two Fano manifolds with semistable tangent bundle.

Then the Fano manifold $X\!=\!Y_1\!\times\!Y_2$ has semistable tangent bundle.
\end{proposition}

\begin{corollary}
Let $X$ be isomorphic to a product of a Del-Pezzo surface with $\PP^1$.

Then $T_X$ is semistable.
\end{corollary}
\begin{pf}
By \cite{Fahlaoui} the tangent bundle of a Del-Pezzo surface is a semistable
vector bundle.
\end{pf}

\begin{example}%
Let $X$ be the  blow-up of $\widetilde{\PP^3}(L)$
 with center
two exceptional lines  of the blow-up
$\widetilde{\PP^3}(L)\to\PP^3$. Here $\widetilde{\PP^3}(L)$ is the blow-up
of $\PP^3$ with center a line $L$.
We will show that $X$ has unstable tangent bundle.

Consider the following diagram
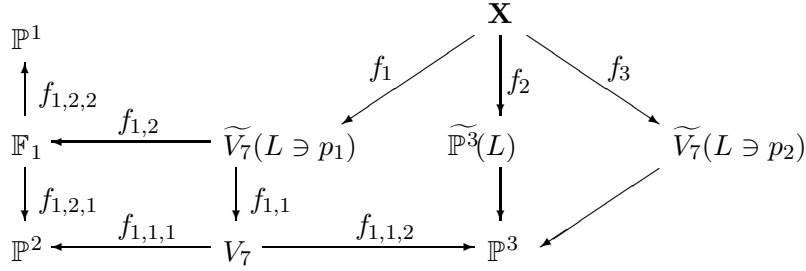
\begin{center}
\begin{picture}(250,100)(50,0)
\put(200,90){$\bf X$}
\put(20,80){$\PP^1$}
\put(20,40){$\FF_1$}
\put(100,40){$\widetilde{V_7}(L\ni p_1)$}
\put(185,40){$\widetilde{\PP^3}\!(L)$}
\put(270,40){$\widetilde{V_7}(L\ni p_2)$}
\put(100,0){$V_7$}
\put(20,0){$\PP^2$}
\put(200,0){$\PP^3$}
\put(195,85){\vector(-3,-2){50}}
\put(155,70){$f\!_1$}
\put(205,85){\vector(0,-1){30}}
\put(207,65){$f\!_2$}
\put(215,85){\vector(3,-2){50}}
\put(245,70){$f\!_3$}
\put(105,35){\vector(0,-1){20}}
\put(110,20){$f\!_{1,1}$}
\put(205,35){\vector(0,-1){20}}
\put(265,35){\vector(-3,-2){45}}
\put(95,45){\vector(-1,0){60}}
\put(60,50){$f\!_{1,2}$}
\put(25,55){\vector(0,1){20}}
\put(30,60){$f\!_{1,2,2}$}
\put(25,35){\vector(0,-1){20}}
\put(30,20){$f\!_{1,2,1}$}
\put(115,5){\vector(1,0){80}}
\put(150,10){$f\!_{1,1,2}$}
\put(95,5){\vector(-1,0){60}}
\put(60,10){$f\!_{1,1,1}$}
\end{picture}
\end{center}
where
\begin{itemize}
\item
$V_7\simeq\PP(\OO_{\PP^2}\oplus\OO_{\PP^2}(1))$, the map $f\!_{1,1,2}$ is the
blow-up
of $\PP^3$ in a point $p_1\in L$ and $f\!_{1,1,1}$ is the canonical projection
$V_7\to\PP^2$.
\item
$f\!_{1,1}:\widetilde{V_7}(L\ni p_1)\to V_7$ is the blow-up of $V_7$
with center the strict transform of $L$ in the blow-up of $V_7\to\PP^3$. The
projection
$f\!_{1,2}:\widetilde{V_7}(L\ni p_1)\to\FF_1$ is a $\PP^1$-bundle.
\item
$f_1:X\to\widetilde{V_7}(L\ni p_1)$ is the blow-up of
$\widetilde{V_7}(L\ni p_1)$ with center an exceptional line of the
blow-up $\widetilde{V_7}(L\ni p_1)\to V_7$.
\end{itemize}
Let $H_1\!=\!f_1^*f_{1,1}^*f_{1,1,1}^*\OO_{\PP^2}(1),$
$H_2\!=\!f_1^*f_{1,1}^*f_{1,1,2}^*\OO_{\PP^3}(1),$
$H_3\!=\!f_1^*f_{1,2}^*f_{1,2,2}^*\OO_{\PP^1}(1)$ and
$D_{f_{1,1}}$, $D_{f_1}$ the pull-backs of the exceptional divisors
of the blow-ups $f_{1,1}$ resp. $f_1$ on $X$. Then we have
\begin{eqnarray*}
&& -K_X=2H_1+2H_2-D_{f_{1,1}}-D_{f_1}, \ H_3\sim H_1-D_{f_{1,1}}\\
&&\mbox{\hspace{-20pt}}
(a_1H_1\!+\!a_2H_2\!-\!a_3D_{f_{1,1}}\!-\!a_4D_{f_1}).(-\!K_X)^2=
 12a_1+15a_2-5a_3-3a_4,\ (-K_X)^3=46.
\end{eqnarray*}
Let $g:=f_{1,2,2}\circ f_{1,2} \circ f_1$.  Since
$g^*\Omega^1_{\PP^1}\simeq \OO_X(-2H_1+2D_{f_{1,1}})$, it follows
that $g^*\Omega^1_{\PP^1} \subset \Omega^1_X$ is a
  $(-K_X)$- destabilising subsheaf, with
$\mu(g^*\Omega^1_{\PP^1})\!=\!-14 > -\frac{46}{3}\!=\!\mu(\Omega^1_X).$
\qed
\end{example}

\noindent
Before we go to the next example we shall collect some auxiliary lemmas.

\begin{Lemma}\label{Conic-bundel.1}
Let $S$ be a smooth projective surface, $f: X \longrightarrow S$
a conic bundle and $\Delta \subset S$ the discriminant locus. Then
we have an exact sequence
$$ 0 \longrightarrow f^*\Omega^1_S \stackrel{\delta}{\longrightarrow}
\Omega^1_X \longrightarrow \Omega^1_{X/S} \longrightarrow 0 $$
and
$\Omega^1_{X/S}\simeq{\cal I}_\Gamma\otimes\omega_X\otimes f^*\omega^{-1}_S$,
where $\Gamma$ is a closed Cohen-Macaulay subscheme of $X$ of pure dimension 1
with $f(\Gamma)=\Delta$. The restriction
$f|_{\Gamma\!\setminus\! f^{-1}(\Delta_{sing})}:
\Gamma\!\setminus\! f^{-1}(\Delta_{sing})\longrightarrow \Delta_{reg}$
is an isomorphism and $\Gamma \cap X_s=(X_s)_{red}$ for all
$s\in \Delta_{sing}$.
\end{Lemma}
\begin{pf}
 $f^*\Omega^1_S \longrightarrow \Omega^1_X$ drops rank in codimension 2,
whence the first three assertions follow from the theory of the
Eagon-Northcott complex \cite{Eagon-Northcott}.
It is also clear that $f(\Gamma)=\Delta$.

$S$ can be covered by affine open sets $U_\alpha$ such that
$f^{-1}(U_\alpha)$ is isomorphic over $U_\alpha$ to the closed subsheme
of $U_\alpha\times\PP^2$ given by a quadratic equation:
$$ g_\alpha:=\sum_{0\leq i \leq j \leq 2} A_{ij}X_iX_j ,\ \ \
A_{ij} \in H^0(U_\alpha,\OO_{U_\alpha}).$$
Using the diagram
$$\begin{array}{ccccccc}
                   &           & & 0 & &  \\
                   &           & & \downarrow & &  \\
0  \longrightarrow & \OO_{U_\alpha\times\PP^2}(-2)|_{f^{-1}(U_\alpha)} &
\longrightarrow & \Omega^1_{U_\alpha\times\PP^2/U_\alpha}|_{f^{-1}(U_\alpha)} &
\longrightarrow & \Omega^1_{X/S}|_{f^{-1}(U_\alpha)} & \longrightarrow  0 \\
&  \multicolumn{2}{r}{(\frac{\partial g_\alpha}{\partial X_0},
\frac{\partial g_\alpha}{\partial X_1},\frac{\partial g_\alpha}{\partial
X_2})\searrow}&
\downarrow  &  \\
&             & &3\OO_{U_\alpha\times\PP^2}(-1)|_{f^{-1}(U_\alpha)} & & & \\
&  &                   & \multicolumn{2}{r}{\downarrow\ \ \ (X_0,X_1,X_2)}  &
\\
&             & & \OO_{U_\alpha\times\PP^2}|_{f^{-1}(U_\alpha)} & &  \\
&                            & & \downarrow & &  \\
&                            & & 0 & &  \\
	\end{array}$$
one can deduce that $\Gamma\cap f^{-1}(U_\alpha)$ is the closed subscheme of
$U_\alpha\times\PP^2$ given by the equations:
$$\frac{\partial g_\alpha}{\partial X_0}=\frac{\partial g_\alpha}{\partial
X_1}=
\frac{\partial g_\alpha}{\partial X_2}=0 $$
(in fact, the  three equations are enough by Euler's identity).
Now the last 2 assertions are clear.
\end{pf}

\begin{Lemma}\label{Cbundel.Aufblasung.T}
Let $Y \stackrel{\pi}{\longrightarrow} X$ be the blow-up of a conic bundle
$ X \stackrel{f}{\longrightarrow} S$ with center a smooth irreducible
subsection $C$
over $S$ (i.e.~ $f|_C:C\to S$ is an embedding). Let $L \in Pic(Y)$, such that
$(\pi_*L)^{**}$ is a
$f$-ample line bundle on $X$.
Then:
\begin{eqnarray*}
&(i)& H^0(Y,T_Y \otimes L^{-1})=0,
 \mbox{ if } \mu(L) > \mu(T_{Y/S}) \mbox{ and}\\
&(ii)& H^0(Y,\Omega^1_Y \otimes L^{-1})=0.
\end{eqnarray*}
\end{Lemma}
\begin{pf}
Straightforward and left to the reader.
\end{pf}

\begin{example}
Let $X$ be the blow-up of $\PP^3$ with center a union of a cubic $C$ in a
plane $S$ and a point $p$ not in $S$.
We will prove that $X$ has stable tangent bundle. For this we make use of the
following diagram\\
\begin{picture}(300,150)(0,0)
\put(145,70){$\bf X$}
\put(0,70){$\PP(\OO_{\PP^2} \oplus \OO_{\PP^2}(1))$}
\put(110,0){$\PP(\OO_{\PP^2} \oplus \OO_{\PP^2}(2))$}
\put(260,70){$Y_3$}
\put(145,130){$Y_4$}
\put(42,130){$\PP^3$}
\put(42,0){$\PP^2 $}
\put(260,0){$W_4$}
\put(140,73){\vector(-1,0){45}}
\put(115,79){$f_1$}
\put(150,65){\vector(0,-1){55}}
\put(154,35){$f_2$}
\put(157,73){\vector(1,0){100}}
\put(198,79){$f_3$}
\put(150,80){\vector(0,1){45}}
\put(154,100){$f_4$}
\put(50,65){\vector(0,-1){55}}
\put(54,35){$\alpha$}
\put(105,5){\vector(-1,0){45}}
\put(85,9){$\beta$}
\put(195,5){\vector(1,0){60}}
\put(210,10){$\gamma$}
\put(263,65){\vector(0,-1){55}}
\put(267,35){$\delta$}
\put(140,133){\vector(-1,0){81}}
\put(95,137){$\varepsilon$}
\put(50,83){\vector(0,1){45}}
\put(54,100){$\zeta$}
\end{picture}\\
\mbox{}

Define $\zeta:V_7=\PP(\OO_{\PP^2}\oplus\OO_{\PP^2}(1))\to \PP^3$ to be
the blow-up of $\PP^3$ in $p$ and denote again by $C\subset V_7$ the
proper transform of $C\subset S \subset \PP^3$. Define also $f_2$ to
by the elementary transformation of $f_1$ along C, and $W_4$ to be
the cone over the Veronese surface in $\PP^5$ and $\gamma$ the
blow-up of the vertex.

Let  $H_1:=f_1^*\alpha^*\OO_{\PP^2}(1)$,
$H_2:=f_1^*\zeta^*\OO_{\PP^3}(1)$,
$H_3:=f_2^*\OO_{\PP(\OO_{\PP^2}\!\oplus\OO_{\PP^2}(2))}(1)$
 and $D_{f_1}$, $D_{f_2} $
the exceptional divisors of $f_1$ resp. $f_2$.
Then it follows that we have
\begin{eqnarray*}
&& D_{f_2} \sim 3H_1 -D_{f_1},\ \ H_2-H_1 \sim H_3-D_{f_2} \\
&& -K_X\sim 2H_1+2H_2-D_{f_1}\sim H_1+2H_3-D_{f_2}\\
&& \Rightarrow\ -K_X\sim H_2+H_3\\
&&(a_1H_1+a_2H_2+a_3H_3).(-K_X)^2= 9a_1+13a_2+19a_3,\ (-K_X)^3=32.\\
&& a_1H_1+a_2H_2+a_3H_3\sim\\
&& \sim (a_1+2a_3)H_1+(a_2+a_3)H_2-a_3D_{f_1}\\
&& \sim (a_1+a_2)H_1+(a_2+a_3)H_3-a_2D_{f_2}.
\end{eqnarray*}
Since $f_2:X\longrightarrow\PP(\OO\oplus\OO(2))$ is the elementary
transformation of
$f_1:X\longrightarrow\PP(\OO\oplus\OO(1))$ the blow-up of
$\PP(\OO\oplus\OO(1))$
with center a smooth subsection $C$ over $\PP^2$, we have an
isomorphism
$D_{f_2}\stackrel{f_1}{\simeq}\alpha^{-1}(C^\prime)\simeq
 Z:=\PP(\OO_{C^\prime}\oplus\OO_{C^\prime}(1))
\stackrel{\rho}{\longrightarrow} C^\prime :=\alpha(C).$
By this isomorphism, the curve $C\subset\alpha^{-1}(C^\prime)$ corresponds to
a curve $C^{\prime\prime}$ on $D_{f_2}$.
The map $g:=\alpha\circ f_1:X\longrightarrow \PP^2$ is a conic bundle with
discriminate locus $C^\prime$.

Since $C\subset Z$ is linear equivalent
to $H_2|Z$ and
$D_{f_1}|D_{f_2}=C^{\prime\prime}\stackrel{f_1}{\simeq} C$, it follows
that
$\OO_X(D_{f_1})|D_{f_2}\simeq\OO_{D_{f_2}}(C^{\prime\prime})\simeq
\OO_Z(C)=\OO_Z(1)$.
Furthermore, we have:
$\OO_X(H_1)|D_{f_2} \simeq \rho^*\OO_{C^\prime}(1)$  and
$\OO_X(H_2)|D_{f_2} \simeq \OO_Z(1)$. Thus, we have
 $\OO_X(H_3)|D_{f_2}\simeq\rho^*\OO_{C^\prime}(2)$ and
$\OO_X(a_1H_1+a_2H_2+a_3H_3)|Z\simeq
\rho^*\OO_{C^\prime}(a_1+2a_3)\otimes\OO_Z(a_2).$

Let $\CL=\OO_X(a_1H_1+a_2H_2+a_3H_3) \subset T_X$.
Consider the exact sequence
$$0\longrightarrow T_{X/\PP^2}\longrightarrow T_X
\longrightarrow g^*T_{\PP^2}.$$
Using lemma
 \ref{Cbundel.Aufblasung.T}, it follows that $a_2+a_3\leq 0$, if
$9a_1+13a_2+19a_3 > 5=\mu(T_{X/\PP^2}).$\\
Let
 $\CL \subset g^*T_{\PP^2}$ be a line bundle with maximal $\mu(\CL)$.
Using
\begin{eqnarray*}
&&  0\not=H^0(X,g^*T_{\PP^2} \otimes \CL^{-1}) \subset
H^0(D_{f_2},g^*T_{\PP^2}\otimes \CL^{-1}|D_{f_2})\simeq\\
&&\simeq
H^0(C^\prime,T_{\PP^2}|_{C^\prime}\otimes\OO_{C^\prime}(-a_1-2a_3)
\otimes S^{-a_2}(\OO_{C^\prime}\oplus\OO_{C^\prime}(1))),\mbox{ it follows} \\
&& a_2\leq 0 \mbox{ and } a_1+a_2+2a_3\leq 1, \mbox{ because }
\deg C^\prime=3.
\end{eqnarray*}
Hence
$\mu(\CL)=9a_1+13a_2+19a_3=9(a_1+a_2+2a_3)+(a_2+a_3)+3a_2\leq 9$.

Now let $\CL \subset \Omega^1_X$.
By lemma \ref{Conic-bundel.1} we have an exact sequence:
\begin{eqnarray*}
&& \mbox{\hspace{-30pt}}
0\longrightarrow g^*\Omega^1_{\PP^2} \longrightarrow \Omega^1_X
\longrightarrow {\cal I}_{C^{\prime\prime}}\otimes\omega_{X/\PP^2}
 \longrightarrow 0\
(\mbox{ with } \omega_{X/\PP^2}\simeq\OO_X(3H_1-H_2-H_3)\ ).\\
&& \mbox{If }H^0(X,\Omega^1_X\otimes\CL^{-1})\not=0 \mbox{ then }
\left\{ \begin{array}{cl}
	 (1) & H^0(X,g^*\Omega^1_{\PP^2}\otimes\CL^{-1})\not=0\mbox{ or}\\
	 (2) & H^0(X,{\cal I}_{C^{\prime\prime}}
	 \otimes\omega_{X/\PP^2}\otimes\CL^{-1})\not=0.
	 \end{array} \right.
\end{eqnarray*}
(1) $g_*(g^*\Omega^1_{\PP^2}\otimes\CL^{-1})\subset
\Omega^1_{\PP^2}(-a_1-2a_3)\otimes
S^{-a_2-a_3}(\OO_{C^\prime}\oplus\OO_{C^\prime}(1))$
$$\Rightarrow a_2+a_3\leq 0\mbox{ and }a_1+a_2+3a_3\leq -2;$$
$\beta_*f_{2*}(g^*\Omega^1_{\PP^2}\otimes\CL^{-1})\subset
\Omega^1_{\PP^2}(-a_1-a_2)\otimes
S^{-a_2-a_3}(\OO_{C^\prime}\oplus\OO_{C^\prime}(2))$
$$\Rightarrow a_2+a_3\leq 0 \mbox{ and }a_1+3a_2+2a_3\leq -2.$$
If $\CL\subset g^*\Omega^1_{\PP^2}$ has maximal $\mu(\CL)$, then we have:
\begin{eqnarray*}
&& 0\not=H^0(X,g^*\Omega^1_{\PP^2}\otimes\CL^{-1})
\subset  H^0(D_{f_2},g^*\Omega^1_{\PP^2}\otimes\CL^{-1}|D_{f_2})\simeq\\
&& H^0(C^\prime,\Omega^1_{\PP^2}|_{C^\prime}\otimes\OO_{C^\prime}(-a_1-2a_3)
\otimes S^{-a_2}(\OO_{C^\prime}\oplus\OO_{C^\prime}(1)))\\
&& \Rightarrow\ a_2\leq 0 \mbox{ and } a_1+a_2+2a_3\leq -2.
\end{eqnarray*}
(2) ${\cal I}_{C^{\prime\prime}}\otimes\omega_{X/\PP^2}\otimes\CL^{-1}\subset
\OO_X((3-a_1)H_1+(-1-a_2)H_2+(-1-a_3)H_3)$
$$\Rightarrow \left\{ \begin{array}{clcl}
		 (\alpha\circ f_1)_*: & a_2+a_3\leq -2
                 & \mbox{and} & a_1+a_2+3a_3\leq -1;\\
		 (\beta\circ f_2)_* : & a_2+a_3\leq -2
                 & \mbox{and} &  a_1+4a_2+2a_3\leq -2.
			 \end{array} \right. $$

If  $0\not=H^0(X,{\cal I}_{C^{\prime\prime}}\otimes
\OO_X((3-a_1)H_1+(-1-a_2)H_2+(-1-a_3)H_3))$, then
$\omega_{X/\PP^2}\otimes\CL^{-1}$ has a global section vanishing
on $C^{\prime\prime}$. It follows that
$\omega_{X/\PP^2}\otimes\CL^{-1}|D_{f_2}$
has a global section vanishing on $C^{\prime\prime}$. Therefore
\begin{eqnarray*}
&&
0\not=H^0(D_{f_2},\OO_X((3-a_1)H_1+(-1-a_2)H_2+(-1-a_3)H_3)
\otimes\OO_{D_{f_2}}(-C^{\prime\prime}))\\
&& \simeq
H^0(Z,\rho^*\OO_{C^\prime}(1-a_1-2a_3)\otimes\OO_Z(-2-a_2)),
 \mbox{ implies that}\\
&& a_2\leq -2 \mbox{ and } a_1+a_2+2a_3\leq -1.
\end{eqnarray*}
Hence
$\mu(\CL)=9a_1+13a_2+19a_3=9(a_1+a_2+2a_3)+(a_2+a_3)+3a_2\leq -17$.
\qed
\end{example}

\end{document}